\documentclass[12pt]{article}

\usepackage[latin1]{inputenc}
\usepackage[T1]{fontenc}
\usepackage[english]{babel}
\usepackage{amsfonts}
\usepackage{euscript}

\newtheorem{conj}{Conjecture}[section]
\newtheorem{theo}[conj]{Theorem}
\newtheorem{prop}[conj]{Proposition}
\newtheorem{coro}[conj]{Corollary}
\newtheorem{lemm}{Lemma}[section]

\newtheorem{defi}{Definition}[section]

\begin{document}

\newcommand{\as}{${\cal A}_6\ $}
\newcommand{\eps}{prehomogeneous symmetric space }
\newcommand{\epss}{prehomogeneous symmetric spaces }
\newcommand{\mt}{${\cal M}_3\ $}
\newcommand{\p}{{\bf I \hspace{-3pt} P}}
\newcommand{\R}{\mathbb{R}}
\newcommand{\C}{\mathbb{C}}
\newcommand{\hh}{\mathbb{H}}
\newcommand{\oo}{\mathbb{O}}
\newcommand{\n}{\mathfrak n}
\newcommand{\B}{\mathfrak b}
\newcommand{\nd}{$\frac{n}{2}\,$}
\newcommand{\nq}{$\frac{n}{4}\,$}
\newcommand{\f}{${\cal F}(Q_1)\ $}
\newcommand{\fp}{${\cal F}'(Q_1)\ $}
\newcommand{\fx}{${\cal F}_x(Q_1)\ $}
\newcommand{\fpx}{${\cal F}_x'(Q_1)\ $}
\newcommand{\se}{Se^k X-X}
\newcommand{\ssi}{if and only if }
\newcommand{\tr}{{}^t}
\newcommand{\trace}{\mbox{tr}}
\newcommand{\im}{\mbox {Im }}
\newcommand{\var}{variety }
\newcommand{\g}{\mathfrak g}
\newcommand{\h}{\mathfrak h}

\newcommand{\dem}{\underline {\bf Proof :} }
\newcommand{\rem}{\underline {\bf Remark :} }
\newcommand{\fin}{\begin{flushright}  $\bullet$ \end{flushright}}
\newcommand{\lpara}{\vspace{-5pt} \ \\}
\newcommand{\para}{\vspace{1pt} \ \\}
\newcommand{\Para}{\vspace{15pt} \ \\}

\newcommand{\fonction}[5]{
\begin{array}{rrcll}
#1 & : & #2 & \rightarrow & #3 \\
   &   & #4 & \mapsto     & #5
\end{array}  }

\newcommand{\fonc}[3]{
#1 : #2 \mapsto #3  }

\newcommand{\matdd}[4]{
\left (
\begin{array}{cc}
#1 & #2  \\
#3 & #4
\end{array}
\right )   }

\newcommand{\mattt}[9]{
\left (
\begin{array}{ccc}
  #1 & #2 & #3 \\
  #4 & #5 & #6 \\
  #7 & #8 & #9
\end{array}
\right )   }

\newcommand{\sectionplus}[1]{\section{#1} \indent \vspace{-5mm}}

\title{Scorza varieties and Jordan algebras}
\author{Pierre-Emmanuel Chaput\\
        chaput@ujf-grenoble.fr}
\date{April 2002}
\maketitle

{\def\thefootnote{\relax}
\footnote{\hskip-0.6cm
{\it AMS mathematical classification \/}: 14M07,14M17,14E07.}}

\begin{center}
\bf{Abstract}
\end{center}

In his book \cite{zak}, F.L. Zak defines and classifies Scorza varieties. A
$k$-Scorza variety is by definition an irreducible smooth 
complex projective variety, of
maximal dimension among those whose $(k-1)$-secant variety is not all of the
ambient space (a precise definition will be given in the second section). 
Let me also recall that the Jordan algebras are the commutative but not
necessarily associative algebras in which the relation\\ $A*(B*A^2)=(A*B)*A^2$
holds. A classical theorem (cf \cite{faraut} for instance) tells us that the
simple real Jordan algebras of rank greater or equal
to 3 (the rank is the generic dimension of
the subalgebra generated by one element) are the algebras of Hermitian 
matrices, with entries in a real normed algebra ($\R,\C,\hh$ or $\oo$), the
product being defined by $A*B=1/2(AB+BA)$ (if the rank is greater or equal to
4, this
normed algebra cannot be the octonions; otherwise the algebra is not a Jordan
algebra). The complex Jordan algebras are the complexifications of real Jordan 
algebras. From Zak's classification theorem we see that there is a very
strong link between Scorza varieties and Jordan algebras:

\begin{theo}[Zak]{The $k$-Scorza varieties are the projectivizations of the 
varieties of rank 1 matrices in the simple Jordan algebras of rank $k+1$.}
\end{theo}

However, this link, as well as the fact that Scorza varieties are homogeneous,
meaning that their automorphism group acts transitively, is not explained by
the geometric proof of F.L. Zak. In \cite{prepub}, I have studied a particular
case of this theorem, namely the 2-Scorza varieties which are also called
Severi varieties, and which have been studied quite a lot because they are the
limiting case of a conjecture by Hartshorne, proved by F.L. Zak. I have shown
how it is possible to adapt Zak's proof so as to get the homogeneity, and to
conclude the proof more easily. The link with Jordan algebra  was nevertheless
not explained in this paper. Here, I wish to give two simple explanations of 
this link, using firstly an old result of McCrimmon\cite{mccrimmon} 
and secondly facts concerning \epss \cite{bertram}.

I wish to thank warmly L. Manivel for helping me on this subject, as well as
W. Bertram and F.L. Zak for precise explanations.

\sectionplus{Scorza and Severi varieties}

In this section, I recall definitions and results from Zak's book. Let
$X\subset \p^N$ be a complex projective variety. Let's denote by $S^kX$ the
closure of the union of all $\p^k$'s containing $k+1$ linearly 
independant points
of $X$. With this notation we have $S^0X=X$ and $S^1X=Sec(X)$, the secant of 
$X$. If $X$ is non degenerate, that is included in no hyperplane, we will call
$k_0(X)$ the least integer $k$ such that $S^kX=\p^N$. In the sequel, all
varieties will be supposed to be non degenerate.
\para

When $X$ and $Y$ have dimension respectively $n$ and $p$, we expect the join
between these two varieties (the closure of the union of lines through a point
of $X$ and a point of $Y$), denoted by $S(X,Y)$, to have dimension $n+p+1$.
Let's denote by $\delta (X,Y)$ and call 'defect of $X$ and $Y$' the difference
between this expected dimension and the true dimension:
$\delta(X,Y):=n+p+1-\dim S(X,Y)$. F.L. Zak introduces for an arbitrary variety
$X$ the defect of $X$ and its secants
$\delta_i:= \delta(X,S^{i-1}X)=\dim X + \dim S^{i-1}X +1 - \dim S^iX$ and calls
$\delta:=\delta_1$. An important result about these defects is theorem 1.8 of
his fifth chapter (\cite{zak} p.109):
\begin{theo}[Zak]{If $1 \leq i \leq k_0$, then
$\delta_i \geq \delta_{i-1}+\delta$.}
\end{theo}

\begin{coro}[Zak]{$k_0 \leq \frac{n}{\delta}$.}
\end{coro}

F.L. Zak then defines Scorza varieties to be the limiting cases of these 
theorems:
\begin{defi}[Scorza varieties]{A variety $X^n \subset \p^N$ is Scorza
if it is irreducible smooth non degenerate and that:
\begin{itemize}
\item{$k_0=\left[ \frac{n}{\delta} \right ]$.}
\item{$\delta_i=i\delta$.}
\end{itemize}  }
\end{defi}
\rem This definition involves the unusual numbers $\delta_i$. But F.L. Zak
shows that if $n$ is a multiple of $\delta$, then the conditions 
$\delta_i=i\delta$ are automatically fullfilled as soon as the first one,
equivalent to the fact that 
$S^{\left [ \frac{n}{\delta} \right ]-1} X\not = \p^N$, is. On the other hand,
during his proof of the classification, he explains that his theorem in the 
case $\delta$ dividing $n$ implies that there is no other Scorza varieties. In
this article, I will suppose that this condition holds.
\para
I call $k$-Scorza variety any Scorza variety such that $k_0=k$. For any $l<k$
and \\
$P\in S^lX-S^{l-1}X$, let $L_P$ denote the closure of the set of points $Q$
in $S^lX$ such that $T_QS^lX=T_PS^lX$ and $Q_P$ the closure of the points
$Q$ in $X$ such that the line $(PQ)$ cuts $S^{l-1}X$ at some point 
different from $Q$.\\
The first step of the proof of the classification theorem is the fact that a
$k$-Scorza variety is made of $(k-1)$-Scorza varieties:
\begin{theo}[Zak]{Let $X^n\subset \p^N$ be any $k$-Scorza variety. Let
$l$ be an integer such that $1\leq l \leq k-1$ and $P$ a generic point 
in $S^lX$.
Then $L_P$ is linear, $Q_P \subset L_P$, and
$Q_P$ is an $l$-Scorza \var in $L_P$. Finally, 
$\dim Q_P = l\delta$ and
$\dim L_P= l(1+\frac{\delta(l+1)}{2})$. \label{zak}}
\end{theo}
\rem I am going to explain why Scorza varieties live in Jordan algebras. In
these algebras, the numbers $k$ and $\delta$ have a nice interpretation, as
well as the previous theorem. In this case $\p^N$ is the projectivised space of
complexified $(k+1)*(k+1)$ Hermitian matrices with entries in a normed algebra.
The dimension of this normed algebra is $\delta$ (hence $\delta$ may only
assume the values 1,2,4 or 8). If $P$ is the usual rank $(l+1)$ matrix
($P=\matdd{I_{l+1}}{0}{0}{0}$), then $L_P$ 
is the linear space of matrices of the 
form $\matdd{*}{0}{0}{0}$, with an arbitrary bloc of dimension $(l+1)*(l+1)$,
and $Q_P$ is the subset of rank 1 matrices. We thus see that $Q_P$ is 
a $l$-Scorza variety.
\para

Lastly, in the sequel, I will use the following two lemmas, which hold in any
$k$-Scorza variety $X$:
\begin{lemm}[\cite{zak},p.127]{$S^{k-1}X$ cannot be a cone. \label{cone}}
\end{lemm}
\begin{lemm}[\cite{zak},p.136]{$S^{k-1}X$ is a hypersurface of degree $k+1$. 
\label{degre}}
\end{lemm}

\sectionplus{Identifying the Jordan algebra}

This section explains the link between Scorza varieties and Jordan algebras
using McCrimmon's theorem 1.2 p.937 in \cite{mccrimmon}. To state this result,
I will need to introduce some notations: if $Q$ is {\it any} 
non-vanishing homogeneous form of
degree $q$ on a linear space $V$, let $G$ be the rational morphism defined by
$$\fonction{G}{V}{V^*}{M}{Q(M,\ldots,M,.)/Q(M)}$$
($Q(A_1,\ldots,A_q)$ denotes the polarisation of $Q$, 
that is the only $q$-linear
symmetric form such that $Q(M,\cdots,M)=Q(M)$). Let also denote by $\tau_M$,
for any $M\in V$ the linear map $-D_MG$ (that is minus the differential of $G$
at the point $M$). For $I$ fixed such that $Q(I)=1$,
$\tau_I$ induces a bilinear symmetric form on $V$, and we denote
$$
\begin{array}{rcl}
\langle A,B \rangle &=& -D_IG(A).B=-D_IG(B).A\\
                    &=& qQ(I,\ldots,I,A)Q(I,\ldots,I,B)-(q-1)Q(I,\ldots,I,A,B)
\end{array}
$$
We may also see that the expression $D_IG(A).B$ is symmetric in $A$ and $B$ by
noting that it is the second differential, evaluated in $I$ on $A$ and $B$, of
$\frac{1}{q}\log[Q(M)]$. If $\tau_I$ is non-degenerate, we define for all $A$
the linear map $H_A$ by the relation  $\tau_A(B,C)=\tau_I(H_AB,C)$.
\para
When $V$ is a Jordan algebra, $I$ its identity and $Q$ its determinant, it is
easily seen that $\tau_I$ in non-degenerate and that $\langle A,B \rangle$
equals $\mbox{tr}(A * B)$. On the other hand, via this duality, $G(M)$
identifies with the inverse of $M$. Finally, $H_AB=-D_AG(B)$ identifies with
$A^{-1}BA^{-1}$, so that $Q[H_A(B)]=Q(A)^{-2}Q(B)$. The following result gives
the reverse implication:

\begin{theo}[McCrimmon]{Let $Q$ be a $q$-form on a linear space $V$ and $I$ a
point in $V$ such that $Q(I)=1$, $\tau_I$ is non-degenerate and 
$Q(H_AB)=h(A)Q(B)$ for
a rational function $h$, as soon as the two members are defined. Then $V$ can
be endowed with a Jordan algebra structure where $I$ is the identity by setting
 $A * B=-{\frac{\partial}{2\partial A}}_{|I}H_A(B)$.
\label{thmccrimmon}}
\end{theo}

Moreover, he shows that this algebra is semi-simple.
Here is the definition of this property: in an algebra, we define the trace of
an element to be the trace of the multiplication by this element. Let
$\trace (M)$ denote this number. This linear form yields a bilinear one by the
formula $(A,B)=\trace (A*B)$. An algebra is said to be semi-simple if this
bilinear form is non-degenerate. 
In our case, it will also follow from a very simple argument 
that our algebra is semi-simple.

We are going to apply theorem \ref{thmccrimmon} to $V=\C^{N+1}$ and $Q$ one
equation of the $(k-1)$-ith secant of our Scorza variety $X$. 
To this end, the main step is to show the following proposition:

\begin{prop}{If $G$ and $\tau$ are defined like previously, then for
$A\not \in S^{k-1}X$, $\tau_A$ is a linear isomorphism between $\p V$ and
$\p V^*$ which maps $X$ onto $Y:=(S^{k-1}X)^*$. \label{isom} }
\end{prop}
In this proposition, I have denoted, for $Z\subset \p V$ a projective variety,
$Z^*\subset \p V^*$ its dual variety, that is the closure for the Zariski
topology of the set of all hyperplanes tangent in one smooth point of $Z$.

\para
\rem A consequence of this proposition is that $Y \simeq X$ is a $k$-Scorza
\var in $\p V^*$. We may already prove that $X$ and $Y$ have the same
dimension: in fact the closure of the preimage of $T_PS^{k-1}X$ for generic
$P$ in $S^{k-1}X$ by the map $P\mapsto T_PS^{k-1}X$ is $L_P$, of dimension
$(k-1)(1+\frac{\delta k}{2})$ by theorem \ref{zak}. Thus the image of this
application, $Y$, must be of dimension
$k(1+\frac{\delta (k+1)}{2})-1-(k-1)(1+\frac{\delta k}{2})=k\delta =n$.

\para
\dem This proof is only a generalisation of the proof I gave in \cite{prepub} 
for the Severi case. Whenever $Z\subset \p V$ is a projective variety, let
$\widehat{Z} \subset V$ denote its cone.

Firstly let $A_1,\ldots,A_k$ be $k$ elements in $\widehat{X}$. Then
$A_1+\ldots+A_k$ is in $\widehat{S^{k-1}X}$, so\\ $Q(A_1+\ldots+A_k)=0$.
Expanding this expression, we deduce that\\$Q(A_1,A_1,A_2,A_3,\ldots,A_k)=0$ 
(cf lemma \ref{degre}). But this formula is linear in $A_i,i\geq 2$, so it
remains true if $A_i,i\geq 2$ are arbitrary elements of $V$ since $X$ is
non-degenerate. This remark will lead to a geometric interpretation of
$\tau_A(x)$, for $x\in X$ and $A\not \in S^{k-1}X$ such that 
$Q(x,A,\ldots,A)\not = 0$. In fact the line $(xA)$ cuts again $S^{k-1}X$
in exactly one point, $x'=A - \frac{Q(A)}{qQ(x,A,\ldots,A)} x$. The tangent
hyperplane at $x'$ to $S^{k-1}X$ is then:
$$Q(A,\ldots,A,.)-\frac{(q-1)Q(A)}{qQ(x,A,\ldots,A)}Q(x,A,\ldots,A,.)=0$$
which is $\tau_A(x)$ up to a constant. So we have proved that
$\tau_A(X)\subset Y$. But these two varieties $X$ and $Y$ have the same 
dimension; since $Y$ is irreducible, to see that $\tau_A$ is onto it
is enough to check that it is generically finite. This is a consequence of the
given geometric interpretation and of the

\begin{lemm}{Let $A\in V$ and $P\in S^{k-1}X$ be generic points (in the sense
of theorem \ref{zak}) such that $A\not \in L_P$. Then $X\cap (L_P+A)-L_P$ is
finite.}
\end{lemm}
\dem
In the linear space $M=L_P+A$, $Q_P$ is a $(k-1)$-Scorza variety 
(theorem \cite{zak}); so let $B$ be the hypersurface of $L_P$ equal to
$S^{k-2}Q_P$ and let ${\cal C}(S^{k-1}X \cap M)$
be the set of irreducible components of $S^{k-1}X \cap M$, then we
have an application:
$$ \fonction{\phi}{(M-L_P) \cap X}{{\cal C}(S^{k-1}X \cap M)}{x}{S(x,B)} $$
Let's recall that $S(x,B)$ is the cone with vertex $x$ and basis $B$. If we
suppose by induction that the classification theorem is proved for the
$(k-1)$-Scorza variety $Q_P$, then $Q_P$ identifies with the 
projective \var of rank 1
matrices and $B$ with that of rank at most $k-1$ matrices. Let me also denote
by $Sing^i B$ the set of matrices with rank at most $k-1-i$, so that the 
singular locus of $Sing^i B$ is $Sing^{i+1} B$ and that $Sing^0 B=B$. We then
get that the application $\phi$ is injective since the singular locus of
$S(x,B)$ is $S(x,Sing^1B)$, so that if $S(x,B)=S(x',B)$ then
$S(x,Sing^i B)=S(x',Sing^i B)$ for all $i$; $i=k-2$ yields 
$S(x,Q_P)=S(x',Q_P)$ and since $Q_P$ is smooth, $x=x'$.
\fin
\rem These isomorphisms $D_AG$ will play an essential role in the sequel. In
a Jordan algebra, we have
$-D_{A^{-1}}G(B)=ABA$. This expression is a classical one in the theory of
Jordan algebra; the linear application $B\mapsto ABA$ is often denoted by 
$P(A)$ and called the quadratic representation. If $M(A)$ denotes the 
endomorphism of the Jordan algebra equal to the multiplication by $A$, $P$
satisfies the identity
$P(A)=2M(A)^2-M(A^2)$ and is involved in numerous important theorems: see for
example theorem 4 p.57 in \cite{koecher} and theorem 2.6 p.47 in
\cite{bertram}.

\begin{coro}{$X$ is homogeneous.}
\end{coro}
\dem Let $A,B\in \p V-S^{k-1}X$. Since $\tau_A$ and $\tau_B$ are
isomorphisms between $X$ and $Y$, $\tau_A^{-1} \circ \tau_B$ is an
automorphism of $X$. Therefore it is enough to check that if $U$ and $W$ are
elements in $X$, there are some $A,B \in \p V-S^{k-1}X$ such that
$\tau_A(U)=\tau_B(W)$. But if $P$
is such that $U$ and $W$ do not belong to $T_PS^{k-1}X$, then it is not
possible that $L_P+U \subset S^{k-1}X$ since in that case $U\in T_PS^{k-1}X$.
Let then $A\in (L_P+U)-L_P-S^{k-1}X$ and $B\in (L_P+W)-L_P-S^{k-1}X$; by the 
previous geometric interpretation, since the lines $(AU)$ and $(BV)$ cut
$S^{k-1}X$ at some point in $L_P$, we have $\tau_A(U)=\tau_B(W)=T_PS^{k-1}X$.
\fin

It is now possible to complete the proof of the classification of Scorza 
varieties in the case where $\delta$ divides $n$: let $I\in V$ such that 
$Q(I)=1$
and $A\in V - \widehat{S^{k-1}X}$. Since $\tau_I$ and $\tau_A$ are isomorphisms
between $S^{k-1}X$ and $S^{k-1}Y$, $H_A=\tau_I^{-1} \circ \tau_A$ is an
automorphism of $S^{k-1}X$, and we can deduce the existence of a rational
function $h$ such that $Q(H_AB)=h(A)Q(B)$. Although this is not necessary; 
we can compute $h$: in fact it is defined out of $\{Q=0\}$, and, when defined,
since $H_A$ is a linear isomorphism, never vanishes. Counting degrees
and taking into account the fact that $h(I)=1$, we can 
deduce that $h(A)=Q(A)^{-2}$. Theorem
\ref{thmccrimmon} then implies that $V$ can be equipped with a Jordan algebra
structure such that $I$ is an identity. 
Since this algebra is semi-simple, it is a direct sum of simple ones, and
$Q$ is the product of the determinants on each summand; since $Q$
is irreducible, this algebra is in fact simple. We thus have proved:

\begin{theo}[Zak]{The $k$-Scorza varieties are the projective varieties
of rank 1 matrices in the complexifications of the real Jordan algebras of
$(k+1)*(k+1)$ Hermitian matrices with entries in one of the four real normed 
algebras (if $k>2$ the octonions are not allowed).}
\end{theo}
\rem In this theorem, nothing distinguishes the 'exceptional' Jordan algebra 
$J_3({\mathbb{O}})$ from the others.

\para

I now want to give a proof, different from McCrimmon's and very simple, 
of the fact that the Jordan 
algebra is semi-simple. To this end, let's definitively identify $V$
and its dual towards the isomorphism $\tau_I$ and let's begin computing the
product: first of all, the relation

\begin{equation}
-D_IG(A)=A=(k+1)Q(I,\ldots,I,A)I-kQ(I,\ldots,I,A,.) \label{Q(I,A,.)}
\end{equation}
yields $Q(I,\ldots,I,A,.)=\frac{1}{k}[(k+1)Q(I,\ldots,I,A)I-A]$.\\
Since $D_AG(B)=\frac{kQ(A,\ldots,A,B,.)}{Q(A)}
-\frac{(k+1)Q(A,\ldots,A,B)}{Q(A)^2}Q(A,\ldots,A,.)$, we may then compute that
$$
\begin{array}{rcl}
 A   *   B &=&  \frac{k(k-1)}{2}Q(I,\ldots,I,A,B,.)\\ 
           & &  - \frac{(k+1)k}{2}[ Q(I,\ldots,I,A) Q(I,\ldots,I,B,.) + Q(I,\ldots,I,B) Q(I,\ldots,I,A,.) ]\\
           & &  + [(k+1)^2 Q(I,\ldots,I,A)Q(I,\ldots,I,B) 
                - \frac{(k+1)k}{2} Q(I,\ldots,I,A,B)]I \\
           &=&  \frac{k(k-1}{2} Q(I,\ldots,I,A,B,.) - 
                \frac{k+1}{2} Q(I,\ldots,I,A) [(k+1)Q(I,\ldots,I,B)I-B]\\
           & &  - \frac{k+1}{2} Q(I,\ldots,I,B) [(k+1)Q(I,\ldots,I,A)I-A]\\
           & &  + [(k+1)^2Q(I,\ldots,I,A)Q(I,\ldots,I,B) -
                \frac{(k+1)k}{2} Q(I,\ldots,I,A,B)]I \\
           &=&  \frac{k(k-1)}{2}Q(I,\ldots,I,A,B,.)
                + \frac{k+1}{2}[Q(I,\ldots,I,A)B + Q(I,\ldots,I,B)A]\\
           & &  - \frac{k(k+1)}{2} Q(I,\ldots,I,A,B)I
\end{array}
$$

In particular, $I$ is indeed an identity for this product:
\begin{lemm}{$I*A=A$.\label{identite}}
\end{lemm}

To prove that our algebra is semi-simple, it is enough to show that $(.,.)$
and $\langle .,. \rangle$ are colinear, since we already know that 
$\langle .,. \rangle$ is non-degenerate. This can be managed in two steps:

\begin{lemm}{$\trace (M)=\frac{(k+1)(2+k\delta)}{2} Q(I,\ldots,I,M)$.}
\end{lemm}
\dem Let us consider the rational function $M\mapsto \det (-D_MG)$. This is
defined away from $\{ Q=0 \}$ and by the following proposition  \ref{isom}, it
never vanishes; moreover it has degree $-(k+1)(2+k\delta)$. Thus
$\det (-D_MG)$ and $Q(M)^{-(2+k\delta)}$ are colinear. Since moreover
$-D_IG \simeq Id$, $\det (-D_MG)=Q(M)^{-(2+k\delta)}$. Differentiating 
this equality yields the lemma.

\begin{lemm}
$\langle A,B \rangle =Q(I,\ldots,I,A * B)$.
\label{scal_prod}
\end{lemm}
\dem First notice that since $Q(I,\ldots,I,.)=-D_IG(I)=I$, 
we have the relation\\
$Q[I,\ldots,I,Q(I,\ldots,I,A,B,.)]=Q(I,\ldots,I,A,B)$. Then
$$
\begin{array}{rcl}
Q(I,\ldots,I,A * B) &=&\frac{k(k-1)}{2}Q[I,\ldots,I,Q(I,\ldots,I,A,B,.)]\\
                    & &+\frac{k+1}{2} Q(I,\ldots,I,A)Q(I,\ldots,I,B)\\
                    & &+\frac{k+1}{2} Q(I,\ldots,I,B)Q(I,\ldots,I,A)
                          -\frac{k(k+1)}{2} Q(I,\ldots,I,A,B) \\
                    &=&(k+1) Q(I,\ldots,I,A)Q(I,\ldots,I,B)
- k Q(I,\ldots,I,A,B)
\end{array}
$$
\fin

\Para
From proposition \ref{isom} we can also deduce a geometric construction of
the product of two elements in $X$:
\begin{prop}{Let $A$ and $B$ be two elements of $\widehat{X}$. Then $A*B$ is
the orthogonal projection of 
$[(k+1)^2Q(I,\ldots,I,A)Q(I,\ldots,I,B)- \frac{k(k+1)}{2} 
Q(I,\ldots,I,A,B)]I$ on\\ $T_A \widehat X \cap T_B \widehat X$.}
\end{prop}
For two generic $A,B \in X$, $T_A \widehat X \cap T_B \widehat X$ does not meet
its orthogonal space (this can be checked using the classification), 
so that the map 'orthogonal projection on
 $A,B \in X$, $T_A \widehat X \cap T_B \widehat X$' is defined generically.\\
\rem
For two generic elements $A$ and $B$, and $P$ generic on the line $(AB)$,
we have by Terracini's lemma 
$T_A \widehat X+T_B \widehat X=T_P \widehat{S^1X}$, so that the vector space
$T_A \widehat X \cap T_B \widehat X$ has dimension $\delta$.

\para
\dem This product $A*B$ is by definition colinear with the derivative at $I$ in
the direction $A$ of the function $\tau_A(B)$. Moreover when
$A$ equals $I$, this element is $B$, and it is in any case an element of $X$.
This proves that $A*B$ belongs to the tangent space $T_B \widehat X$. 
By symmetry, it is also in $T_A \widehat X$.

For $U \in T_AX \cap T_BX$, let us compute the scalar product 
$\langle U,A*B \rangle$. As we have seen that the elements $A$ of $X$ satisfy
$Q(A,A,.,\ldots,.)=0$; and since $U \in T_AX \cap T_BX$, then
$Q(A,U,.,\ldots,.)=Q(B,U,.,\ldots,.)=0$. Then\\
$\langle A,U \rangle = (k+1)Q(I,\ldots,I,A)Q(I,\ldots,I,U)$. We deduce:
$$
\begin{array}{rcl}
\langle U,A*B \rangle &=& \langle U, \frac{k(k-1)}{2}Q(I,\ldots,I,A,B,.)\\
     &&         + \frac{k+1}{2}[Q(I,\ldots,I,A)B + Q(I,\ldots,I,B)A]
                - \frac{k(k+1)}{2} Q(I,\ldots,I,A,B)I \rangle\\
                    &=& (k+1)^2Q(I,\ldots,I,A)Q(I,\ldots,I,B)Q(I,\ldots,I,U)\\
     &&         - \frac{k(k+1)}{2} Q(I,\ldots,I,A,B)Q(I,\ldots,I,U)\\
                    &=& \langle U,[(k+1)^2Q(I,\ldots,I,A)Q(I,\ldots,I,B)
                - \frac{k(k+1)}{2} Q(I,\ldots,I,A,B)]I \rangle\\
\end{array}
$$
\fin

For fixed $A\in X$, this proposition only defines $M_A$ on $X$ up to a 
constant, but then it does not define $M_A$ on $V$, even up to a constant. The
following lemma and proposition enable us to give a precise geometric
definition of the product.

\begin{lemm}
Let $V$ and $W$ be vector spaces; $\varphi,\psi:\p V \rightarrow \p W$
linear maps and $X \subset \p V$ a non-degenerate irreducible variety.
Suppose that the projective applications 
$\varphi, \psi : X \rightarrow \p W$ are the same and that
$\ker \widehat \varphi = \ker \widehat \psi$. 
Then $\varphi$ and $\psi$ are projectively equal.
\label{app_lin}
\end{lemm}
$\widehat \varphi$ and $\widehat \psi$ are the linear maps 
$V \rightarrow W$ corresponding to $\varphi$ and $\psi$.\\
\rem Be careful that we make the strong hypothesis that 
$\widehat \varphi$ and $\widehat \psi$
have the same kernel in $V$, and not only in $\widehat X$.\\
\dem We have 
$\im \widehat \varphi = \langle \widehat \varphi(\widehat X) \rangle = 
\langle \widehat \psi(\widehat X) \rangle = \im \widehat \psi=:I$ 
and by hypothesis 
$\ker \widehat \varphi = \ker \widehat \psi=:K$. 
Thus there exists $A \in GL(I)$ such that 
$\widehat \psi = A \widehat \varphi$. 
In fact if $V'$ is a complementary space of $K$ in $V$ and
$\varphi',\psi'$ are the isomorphisms of $V'$ with $I$ such that, on 
$V=K \oplus V'$, we have $\widehat \varphi = (0,\varphi')$ and 
$\widehat \psi=(0,\psi')$, it is
enough to set $A=\psi' \circ {\varphi'}^{-1}$. Now, since $\varphi$ and $\psi$
are equal on $X$, $\psi(X)$ only contains eigenvectors of $A$. But
since $\psi(X)$ is irreducible and non-degenerate in $I$, $A$ has only one
eigenspace, $I$; it is thus a multiple of the identity and $\varphi$ et $\psi$
are projectively equal.
\fin

\begin{prop}
Let $A \in X,B \in \p V$. Then $A*B=0$ \ssi $B\in \Sigma_{A'}$.
\label{produit_nul}
\end{prop}
This proposition and the previous one define, according to the previous lemma,
the algebra structure.\\
\dem 
By proposition \ref{isom}, for $B \not \in S^{k-1}X$, to construct
$\tau_B(A)$, we first need to construct the intersection point $A''$ of the 
line $(AB)$ with $S^{k-1}X$, and then say that $\tau_B(A)$ is the hyperplane
$T_{A''} S^{k-1}X$. Similarly, $\tau_I(A)$ is the hyperplane
$T_{A'} S^{k-1}X$, where $A'$ is the intersection point of the line $(IA)$
with $S^{k-1}X$. Thus, we have $Q(B).A \parallel A$ \ssi these two hyperplanes
are the same, which means that
$\Sigma_{A''}=\Sigma_{A'}$ or
$A''\in \Sigma_{A'}$, or thus that
$B \in \langle \Sigma_{A'},I \rangle$. We thus have proved
$$Q(B).A \parallel A \Longleftrightarrow B \in \langle \Sigma_{A'},I \rangle$$
Taking the tangent spaces to these varieties at $I$, we deduce
$$A*B \parallel A \Longleftrightarrow B \in \langle \Sigma_{A'},I \rangle$$

Suppose now that $\langle A,I \rangle \not = 0$. Then we have
$\langle \Sigma_{A'},I \rangle \cap A^\bot = \Sigma_{A'}$. In fact, these two
sets are hyperplanes in $\langle \Sigma_{A'},I \rangle$ and we always have
$\Sigma_{A'} \subset A^\bot$, since $A^\bot=T_{A'}S^{k-1}X$ by
definition of the scalar product. Moreover, since it is clear that if
$A*B=0$, then $\langle A,B \rangle = \langle A*B , I \rangle = 0$ 
(lemma \ref{scal_prod}), 
we thus have for generic $A$:\\
$\{ B : A*B = 0 \} = \Sigma_{A'}$. For $A$ such that
$\langle A,I \rangle = 0$, $\{ B : A*B = 0 \}$ is again a hyperplane in
$\{ B : A*B \parallel A \} = \langle \Sigma_{A'},I \rangle$, 
by continuity it is still $\Sigma_{A'}$.
\fin

\sectionplus{Identification of $V$ with a \eps}

In the matrix algebras we have met in the previous section, the open subset
of invertible matrices is homogeneous under the action of the group preserving
the determinant up to a constant. Such a vector space equipped with a group 
action such that there exists a dense orbit is called prehomogeneous. If
moreover there is an involution of the group such that the stabilizer of a 
point
in the dense orbit is included in the set of fixed points for this involution
and contains the connected component of the identity of this set, 
then we say that this space is symmetric.
\para

The goal of this section is to prove that if $X\subset \p V$ is a $k$-Scorza
variety, then $V$ is a \eps (without using the classification theorem), which
provides another proof of the classification of Scorza varieties. Finally, as
far as the fourth Severi variety $X\subset \p J_3(\mathbb{O} )$
is concerned, we get that the quasi-projective
variety of invertible matrices identifies with $E_6/F_4$, where $F_4$ is the
subgroup of the adjoint group $E_6$ preserving a non-degenerate quadratic form.
\Para

So let  $X\subset \p V$ be a Scorza variety, and as in the previous section
$Q$ an equation of  $S^{k-1}X$, $G$ defined by 
$G(M)=\frac{Q(M,\ldots,M,.)}{Q(M)}$, $I$ such that $Q(I)=1$ and $*_K$ defined 
by $A*_KB=\frac{\partial}{2\partial A}_{|I}[D_AG(B)]$. Since there will soon
be another product, I want to put an index $*_K$ so as to avoid confusions. Let
finally $\cal G$ be the subgroup of $GL(V)$ preserving $\widehat{S^{k-1}X}$.
\para

\begin{lemm}{$V - \widehat{S^{k-1}X}$ is homogeneous under the action of
$\cal G$.}
\end{lemm}
\dem To prove this result, we consider the elements $\tau_B^{-1} \circ \tau_A$
of $\cal G$. Whenever we are given two elements 
$U,W \in V- \widehat{S^{k-1}X}$, if there exists 
$A,B \in V - \widehat{S^{k-1}X}$ such that $\tau_A(U)=\tau_B(W)$, then $U$ and
$W$ are in the same $\cal G$-orbit. But the set of the 
$\tau_I^{-1} \circ \tau_A(I)$ contains an open set in $V$ since the 
differential of the morphism $A \mapsto \tau_I^{-1} \circ \tau_A(I)$ 
is the application
$A\mapsto -2A*I$, which is $-2Id$ by lemma \ref{identite}, so is invertible.
We thus get that the $\cal G$-orbit of $I$ is dense; since $I$ may have been
chosen to be any element such that $Q(I)=1$, this result holds for any
element $U$ such that $Q(U)\not = 0$, and the proposition follows.
\fin

\begin{prop}{$V$ is a \eps under the action of $\cal G$. 
\label{symetrique}}
\end{prop}
\dem We have just seen that it is prehomogeneous. If $l$ is a linear form, the
equation $Q(l)=0$ is equivalent to $D_IG^{-1}(l) \in S^{k-1}X$, or 
$l \in S^{k-1}Y$. Let us see now that if $g\in \cal G$, then 
$\tr g \in \cal G$: to this end it is enough to see that $\tr g$ preserves
$S^{k-1}Y$. Let $l\in S^{k-1}Y$. For arbitrary $A$ in $\p V-S^{k-1}X$, there
exists $B\in S^{k-1}X$ such that $l=D_AG(B)$. But
$\tr g.D_AG(B)=D_{g^{-1}A}G(g^{-1}B)$, so that $\tr g.l \in S^{k-1}Y$.
\para
So let now define on $\cal G$ the involution 
$\sigma :g\mapsto \tr g^{-1}$. Then $g$ is a fixed point of this involution
\ssi $g$ preserves the quadratic form $\langle A,B \rangle =-D_IG(A)(B)$.\\
We may then finish the proof of proposition \ref{symetrique} by

\begin{prop}{Let $g\in \cal G$. Let us consider the three following conditions:
\begin{enumerate}
\item{$g$ preserves the product on $V$, that is $(g.A)*_K(g.B)=g.(A*_KB)$.}
\item{$g$ preserves $I$, $g.I=I$.}
\item{$g$ preserves the scalar product, 
$\langle g.A,g.B \rangle =\langle A,B \rangle$.}
\end{enumerate}
Then 1 and 2 are equivalent. Moreover, 1 implies 3 and the elements of the
connected component of the identity in the set of elements 
satisfying 3 satisfy 1.}
\end{prop}
{\bf \underline{Remarks :}}
\begin{itemize}
\item{As announced, this proposition identifies the stabilizer $F_4$ of $I$ in
$E_6$ with the group of automorphisms of the exceptional Jordan algebra
$J_3({\mathbb O})$, and also
with the subgroup preserving a non-degenerate quadratic
form.}
\item{The equivalence between 1 and 2 in a Jordan algebra $J$ is expressed by 
the classical fact that an element of the structural group Str($J$) is in
Aut($J$) \ssi $g.I=I$, cf for instance \cite{faraut}, p.148.}
\end{itemize}
\dem
\begin{itemize}
\item{If 1 holds, then $(g.I)*_K(g.B)=g.B$, so that $g.I$ is an identity. But
in an algebra, there can be at most one identity, since if $I$ and $I'$ are
identities, then $I*_KI'=I=I'$. Thus $g.I=I$ (2).}
\item{If 2 holds, since moreover $g$ preserves $Q$ up to a constant, $g$
preserves exactly $Q$. Since
$\langle A,B \rangle=kQ(I,\ldots,I,A,B)-(k+1)Q(I,\ldots,I,A)Q(I,\ldots,I,B)$, 
then
$\langle g.A,g.B \rangle =\langle A,B \rangle$ (3). 
We are going to see that $g$ commutes with $G$ before
deducing 1. If $M$ and $x$ are element of $V$, then\\
$Q(g.M,\ldots,g.M,x)=Q(M,\ldots,M,g^{-1}x)$, so\\
$Q(g.M,\ldots,g.M,.)=\tr g^{-1}.Q(M,\ldots,M,.)=g.Q(M,\ldots,M,.)$. 
Dividing this equality by $Q(M)$, one deduces $G(g.M)=g.G(M)$. Differentiation
yields $D_{g.A}G(g.B)=g.D_AG(B)$ and another differentiation in $I$ gives
$(g.A)*_K(g.B)=g.(A*_KB)$.
\lpara

We can also give another proof, more geometric, of this implication. In fact
let us consider for $g$ satisfying (2) and $A \in \widehat{X}$ fixed
the endomorphisms of $V$ equal to $B \mapsto (g.A)*(g.B)$ and 
$B \mapsto g.(A*B)$. They agree projectively on $X$ since generically they
are the orthogonal projection of $I$ on 
$T_{g.A}X \cap T_{g.B}X=g.(T_AX \cap T_BX)$. Moreover they vanish respectively
when $g.B \in \Sigma_{(g.A)'}$ and when
$B \in \Sigma_{A'}$ (if $x\in X$, then $x'$ denotes the intersection point
between $(xI)$ and $S^{k-1}X$ different from $x$), 
both conditions being equivalent because since $g$ 
preserves $I$ by hypothesis, $g.A'=(g.A)'$. 
By lemma \ref{app_lin}, these two endomorphisms
are proportionnal, and since they agree on $I$ they are equal, and by 
linearity $(g.A)*(g.B)=g.(A*B)$ holds even if $A$ is arbitrary.}
\item{The last thing to be proved is that any element in the component of
the identity in 3 satisfies 2. Any 
$g\in \cal G$ satisfies $Q(g.A)=\lambda Q(A)$
for some scalar $\lambda$. As before, if it satisfies 3, one deduces that
$G(g.M)=g.G(M)$. For $M=I$, this yields $G(g.I)=g.I$. If $g_0\in \cal G$ is
such that $g_0.I=I$, and if $g$ is near $g_0$, since $I$ is an isolated
fixed point of $G$ (in fact the differential in $I$ of the application
$\fonc{\varphi}{N}{G(N)-N}$ is $-3 Id$), then $g.I=I$. We have shown that the
set of all $g$ such that $g.I=I$ is open in the set of $g$ satisfying 3.
\lpara

This implication has again a geometric proof: let $g$ preserve the scalar
product, meaning that for all $A\in \p V$,
the orthogonal hyperplane to $g.A$ is the image by $g$ of the orthogonal
hyperplane to $A$. Choose $A \in X$. Then the former
(respectively the latter) hyperplane is by definition of the scalar 
product the tangent hyperplane to $S^{k-1}X$ 
at the intersection point $(g.A)'$ 
(respectively $g.A'$)
of the line $(g.A,I)$ (respectively $(g.A,g.I)$) 
with $S^{k-1}X$. If they are equal, this implies that 
$(g.A)' \in \Sigma_{g.A'}$, and so that
$g.I \in \langle \Sigma_{(g.A)'},I \rangle$. 

Suppose now that $g.I \not = I$ (projectively) and let $A \in X$. 
Then the line $(I,g.I)$ in $\langle \Sigma_{(g.A)'},I \rangle$ cuts
$ \Sigma_{(g.A)'} $. Thus there exists a point in
$(I,g.I) \cap \Sigma_{(g.A)'} $, thus in
$(I,g.I) \cap T_{(g.A)'}S^{k-1}X \cap S^{k-1}X $. Since 
$(I,g.I) \cap S^{k-1}X$ is finite and since for
fixed $M$ the fact that $M \in T_{(g.A)'}S^{k-1}X $ is a closed condition upon
$A$, there exists a fixed $M$ which lies in every 
$ T_{(g.A)'}S^{k-1}X = \tau_I(g.A)$, 
and thus in all
$T_PS^{k-1}X$ (proposition \ref{isom}), which is contradictory.}
\end{itemize}
\fin

\para
On such a prehomogeneous symmetric space, 
one can define naturally an algebra structure; let us do it
following the notations of W. Bertram \cite{bertram} and see that this product
is nothing else than $*_K$. For the moment, I will denote by $*_B$ this new
product. Let $Id$ be the identity in $\cal G$, then $d\sigma_{Id}$ is an
involution of the Lie algebra $\g$ of $\cal G$; let us
denote by $\g^+$ and $\g^-$
the two eigenspaces associated to the eigenvalues 1 and -1. Let $E$ denote the
evaluation $g\in G \mapsto g.I$ and $e$ its differential. The latter yields an
isomorphism between $\g^-$ and $V$ since the orbit of $I$ is open. The product
on $V$ is then given by $A*_BB=e^{-1}(A).B$. A classical result is that
equipped with this product, $V$ is Lie triple system, meaning that if $A$ and
$B$ are two elements of $V$, then the commutator $[M_A,M_B]$ is a derivation of
$*_B$, recalling that $M_U$ stands for the morphism of multiplication by $U$ 
(cf \cite{bertram}, p.44).
\para
Let us see that $*_B=*_K$. The elements $D_AG$ of $\cal G$ are symmetric in the
sense that $D_AG(B)(C)=D_AG(C)(B)$, as was shown in the preceding section. Thus
they satisfy $\tr g=g$, or $\sigma(g)=g^{-1}$. The derivative $M_A$ of 
$\frac{-1}{2}D_AG$ at $I$ in the direction $A$ then belongs to $\g$ and 
satisfies $d\sigma_{Id}(g)=-g$, so it is in $\g^-$. Since $A*_KI=A$,
$e^{-1}(A)=M_A$ and $$A*_BB=M_A(B)=A*_KB$$
\para

To conclude that this algebra structure is a Jordan algebra, I can use the
following result (\cite{bertram} p.108):

\begin{theo}{A semi-simple Lie triple system is a Jordan algebra.}
\end{theo}

\sectionplus{The Severi case}

In the case of Severi varieties, one can show by hand a Cayley-Hamilton type
relation, which implies that our algebra is power associative, meaning that
the subalgebra generated by one element is associative. Let me recall that in
that case $Q$ has degree 3.

First I prove that the matrix $Q(A,A,.)$ is the comatrix of $A$:

\begin{lemm}{$A * Q(A,A,.) = Q(A) I$.}
\end{lemm}
\dem
It is a known fact (cf for instance \cite{prepub}) that the birational map
$G$ is an involution if we identify $V$ and $V^*$. Moreover this involution
sends $Sec(X)$ onto $X$. So $x \in Sec(X)$ \ssi $Q(x,x,.)\in Sec(X)$. Then
$Q[Q(A,A,.)]=Q(A)^2$. Differentiating this relation twice will lead to 
relations which imply the lemma. A first derivation yields:\\
$Q[Q(A,A,.),Q(A,A,.),Q(A,B,.)]=Q(A)Q(A,A,B)$. If $U$ is the linear form\\
$Q(A,B,.)$, one reads the relation as 
$Q[Q(A,A,.),Q(A,A,.),U]=Q(A)\langle A,U \rangle$. But
for generic $A$, the map $B\mapsto Q(A,B,.)$ is an isomorphism (it is the case
as soon as $A\not \in Sec(X)$: cf \cite{prepub}); we can deduce that for
generic $A$, the previous relation holds for all $U$, so that
\begin{equation}
Q[Q(A,A,.),Q(A,A,.),.]=Q(A)A. \label {rel1}
\end{equation}
Thus this relation is true for any $A$.
Differentiating again yields
\begin{equation}
4Q[Q(A,A,.),Q(A,B,.),.]=3Q(A,A,B)A+Q(A)B. \label {rel2}
\end{equation}
Letting $B=I$ and taking into account (\ref{Q(I,A,.)}) (p.3) we get
\begin{equation}
2Q[Q(A,A,.),A,.]=6Q(I,I,A)Q[Q(A,A,.),I,.]-Q(A)I-3Q(A,A,I)A. \label{rel3}
\end{equation}
Finally, computing this linear form in $I$, one gets
\begin{equation}
2Q[Q(A,A,.),A,I]=3Q(I,I,A)Q(A,A,I)-Q(A). \label{rel4}
\end{equation}

\para
Now, by definition,
$$
\begin{array}{rcl}
\langle A*Q(A,A,.),U \rangle &=& Q[Q(A,A,.),A,U]+\frac{3}{2}Q(A,A,I)\langle 
A,U \rangle\\
               & & + \frac{3}{2}Q(I,I,A)Q(A,A,U)-3Q[Q(A,A,.),A,I]\langle I,U
\rangle
\end{array}
$$
Taking into account relations (\ref{rel3}) and (\ref{rel4}), one proves then 
that
$$
\begin{array}{rcl}
\langle A * Q(A,A,.),U \rangle &=& 3 Q[Q(A,A,.),I,U)Q(I,I,A) + Q(A)
\langle I,U \rangle \\ 
                 & & + \frac{3}{2} Q(I,I,A)Q(A,A,U) - 
\frac{9}{2} Q(A,A,I)Q(I,I,A)\langle U,I \rangle
\end{array}
$$
Moreover,
$$
\begin{array}{rcl}
Q[Q(A,A,.),I,U] &=& Q[Q(A,A,.),I,3Q(I,I,U)I-2Q(I,U,.)]\\
                &=& \frac{3}{2}Q(A,A,I)Q(I,I,U) - \frac{1}{2}Q(A,A,U)
\end{array}
$$
applying relation (\ref{rel2}) with $A=I$ and $B=U$.

Therefore $\langle A * Q(A,A,.),U \rangle = Q(A)\langle I,U \rangle$. 
\fin

Now, since $Q(A-\lambda I)=Q(A) - 3Q(A,A,I) \lambda + 3Q(A,I,I) 
\lambda^2 - \lambda^3$ is the caracteristic polynomial of $A$, one expects
that:
\begin{lemm}[Cayley-Hamilton]{$A * A * A = 3Q(A,I,I)A^2 - 3Q(A,A,I) A + 
Q(A) I$ \label{cayley}}
\end{lemm}
\dem Since $A * A=Q(A,A,.) + 3 Q(I,I,A) A - 3 Q(I,A,A)I$, this is a consequence
of the previous lemma.
\fin

An easy consequence of this proposition is:
\begin{coro}{The product is power associative.}
\end{coro}
\dem If $I,A,A^2$ are independent vectors in $V$, let $M_A$ denote the matrix,
expressed in this base, of the multiplication by $A$ in the subalgebra 
generated by $A$:
$$M_A=\mattt{0}{0}{Q(A)}{1}{0}{-3Q(A,A,I)}{0}{1}{3Q(A,I,I)}$$
If $I,A$ and $A^2$ are dependent let me define $M_A$ by the same formula. Let
also $A_i$ be the vector which can be expressed as the linear combinaison of
$I,A$ and $A^2$ given by the first column-vector of the matrix $M_A^i$. Then it
is sufficient to convince oneself that any expression written only with $A$, 
the sign $*$ and brackets, equals $A_i$, $i$ being the number of written $A$'s.
If $i\leq 2$, this is trivial; for $i=3$, this a consequence of the previous
lemma; for $i=4$, one has to show that
$$
\begin{array}{rcl}
A^2*A^2=A*A^3 &=& [9Q(A,I,I)^2-3Q(A,A,I)]A^2\\
              & & +[Q(A)-9Q(A,I,I)Q(A,A,I)]A+3Q(A,I,I)Q(A)I
\end{array}
$$
and it is an easy consequence of the shown relations.
Lastly, for $i\geq 5$, one can argue by induction.
\para

A consequence of this corollary and of theorem 2.15 p.52 of \cite{bertram},
which states that a unitary Lie triple system is a Jordan algebra \ssi it is
power associative, is again that our algebra is a Jordan algebra. 
Moreover, lemma
\ref{cayley} shows that this algebra is of rank at most 3.

\vfill
\begin{flushleft}
Pierre-Emmanuel CHAPUT\\
INSTITUT FOURIER\\
Laboratoire de Math\'ematiques\\
UMR5582 (UJF-CNRS)\\
BP 74\\
38402 St MARTIN D'HERES Cedex (France)\\
\vspace{.5cm} chaput@ujf-grenoble.fr
\end{flushleft}

\end{document}